\documentclass[12pt,twoside]{amsart}
\usepackage{amsthm,amsfonts,amssymb,amscd,euscript}

\usepackage[matrix,arrow]{xy}

\author{Florin Ambro} 
\address{Department of Mathematical Sciences\\
University of Tokyo,
Komaba, Meguro-Ku,Tokyo 153-8914, JAPAN.}
\email{ambro@ms.u-tokyo.ac.jp}

\address{Current Address: DPMMS, CMS\\
University of Cambridge,
Wilberforce Road, Cambridge CB3 0WB, England.}
\email{f.ambro@dpmms.cam.ac.uk}

\newcommand{\isoto}{{\overset{\sim}{\rightarrow}}}
\newcommand{\bull}{{\scriptscriptstyle{\bullet}}}

\newcommand{\Q}{{\mathbb Q}}
\newcommand{\Z}{{\mathbb Z}}
\newcommand{\N}{{\mathbb N}}
\newcommand{\R}{{\mathbb R}}

\newcommand{\cF}{{\mathcal F}}

\newcommand{\cE}{{\mathcal E}}
\newcommand{\cO}{{\mathcal O}}

\newcommand{\cI}{{\mathcal I}}

\newcommand{\cS}{{\mathcal S}}

\newcommand{\Amp}{\operatorname{Amp}}
\newcommand{\Div}{\operatorname{Div}}

\newcommand{\Supp}{\operatorname{Supp}}

\newcommand{\mult}{\operatorname{mult} }

\newcommand{\Bsl}{\operatorname{Bsl}}

\theoremstyle{plain}
\newtheorem{thm}{Theorem}[section]

\newtheorem{lem}[thm]{Lemma}
\newtheorem{cor}[thm]{Corollary}

\newtheorem{prop}[thm]{Proposition}

\theoremstyle{definition}
\newtheorem{defn}[thm]{Definition}

\newtheorem{exmp}[thm]{Example}

\newtheorem{rem}[thm]{Remark}

\newtheorem{ack}{Acknowledgments}   

\theoremstyle{remark}

\begin{document}

\bibliographystyle{amsalpha+}
\title[quasi-log varieties]
{Quasi-log varieties}

\begin{abstract} We extend the Cone and Contraction 
Theorems of the Log Minimal Model Program to log 
varieties with arbitrary singularities. 
\end{abstract}

\maketitle


\setcounter{section}{-1}

%

\section{Introduction}


 The starting point of the Minimal Model Program is the Cone 
 and Contraction Theorems of S. Mori: the $K_X$-negative 
 part of the cone of effective curves of a non-singular 
 projective $3$-fold $X$ is locally rationally polyhedral, 
 with contractible faces. 
One hopes that by replacing the original variety with the target 
space of the contraction associated to a negative face, or a small 
modification of it (a flip), we reach a minimal model or a 
Mori-Fano fiber space, after finitely many steps.
These intermediate varieties have singularities in dimension at 
least three, so it became clear that one must consider varieties 
with some mild singularities in order to find minimal models. 

In characteristic zero, Y. Kawamata, X. Benveniste, M. Reid, 
V. V. Shokurov and J. Koll\'ar proved the Cone and Contraction 
theorems for varieties with Kawamata log terminal singularities.
This part of the Log Minimal Model Program was expected
to work for log varieties with arbitrary singularities, under
certain assumptions on rays or their contractions. This 
is our main result, and before we state it we make the 
following definition:

\begin{defn} A {\em generalized log variety} $(X,B)$ is 
a pair consisting of a normal variety $X$ and an effective
Weil $\R$-divisor $B$ such that $K+B$ is $\R$-Cartier.
We denote by $(X,B)_{-\infty}$ the locus were $(X,B)$
does not have log canonical singularities (it has a 
natural subscheme structure). 
A {\em log variety} is a generalized log variety which
has log canonical singularities, i.e. $(X,B)_{-\infty}
=\emptyset$.
\end{defn}

\begin{thm} Let $(X,B)$ be a projective generalized log 
variety defined over a field of characteristic zero. Let 
$\overline{NE}(X)$ be the closure of the cone of effective 
curves of $X$, and set 
$$
\overline{NE}(X)_{-\infty}=
Im(\overline{NE}((X,B)_{-\infty}) \to \overline{NE}(X))
$$
\begin{itemize}
\item[(i)] Let $F$ be a face of the cone $\overline{NE}(X)$
 such that 
 $$
 F \cap (\overline{NE}(X)_{-\infty} +
 \overline{NE}(X)_{K+B \ge 0}) =\{0\}
 $$ 
 Then there exists a projective contraction $\varphi_F:X\to Y$ 
 which contracts exactly the curves belonging to $F$. Furthermore,
 $\varphi_F$ restricted to $(X,B)_{-\infty}$ is a closed 
 embedding.
 \item[(ii)]
$
\overline{NE}(X)=\overline{NE}(X)_{K+B \ge 0}+
\overline{NE}(X/S)_{-\infty}+\sum R_j,
$
where the $R_j$'s are the one dimensional faces satisfying
the assumption in (i). Furthemore, the $R_j$'s are discrete 
in the half space $N_1(X)_{K+B<0}$.
\end{itemize}
\end{thm}

This result is a special case of Theorem~\ref{cone}.
As a corollary, we generalize a result of J. Koll\'ar ~\cite{bug}
(in characteristic zero): 
if $(X,B)$ has log canonical singularities outside a finite set of 
points, the Cone Theorem holds exactly as in the Kawamata log 
terminal case. In particular, this holds for a normal surface
with $\Q$-Gorenstein singularities (cf. ~\cite{Sak}).
See also \cite{PLflips} for applications.

We also establish the Base Point Free Theorem for generalized
log varieties, including the log big case (Theorems~\ref{bpf}, 
~\ref{logbig}). Another application is the uniqueness of minimal
lc centers of (quasi-) log Fano varieties (Theorem~\ref{u}).

For the proof, it turns out to be easier to work in a larger class 
of varieties that we call {\em quasi-log varieties}. Their 
definition is motivated by Y. Kawamata's X-method, which produces global 
sections of adjoint line bundles $L$: we first create singularities 
which are not Kawamata log terminal inside $X$, i.e. $LCS(X)\ne 
\emptyset$. By adjunction, we expect that $L|_{LCS(X)}$ is still 
an adjoint line bundle, hence if it has a global section (by induction, 
for instance), we can lift it to a global section of $L$ by the
Kawamata-Viehweg vanishing. 
Unlike the given variety, its LCS locus is no longer normal, 
not even irreducible or equi-dimensional, and its log canonical 
class in the usual sense does not make sense either. However, 
by definition, the LCS locus is the target space of a $0$-log 
contraction (cf. ~\cite[3.27(2)]{PLflips})
from a variety with only embedded normal crossings singularities.
We call {\em quasi-log varieties} 
those varieties appearing as the target space of such
contractions. Examples are varieties with embedded
normal crossings singularities, generalized log varieties and 
their LCS loci (see ~\ref{basicexmp}).

A quasi-log variety $X$ is endowed with an $\R$-Cartier
divisor $\omega$, the descent of the log canonical class of 
the total space of the $0$-log contraction, a closed proper 
subscheme $X_{-\infty}\subset X$, and a finite family $\{C\}$ 
of reduced and irreducible subvarieties of $X$. 
We say that $\omega$ is the {\em quasi-log canonical class} of 
$X$, $X_{-\infty}$ is the locus where $X$ does not have 
{\em qlog canonical} singularities, and the $C$'s are the
{\em qlc centres} of $X$. 
The open subset $X\setminus X_{-\infty}$ is reduced, with 
seminormal singularities. 
We note here that singularities appearing on special 
LCS loci have been called {\em semi-log canonical} in the 
literature.

The adjunction and vanishing for quasi-log varieties 
are proved in Theorem~\ref{adj_van}. The former holds 
by the very definition, while the latter is an extension 
to normal crossings pairs of the vanishing and torsion 
freeness theorems of J. Koll\'ar, based on previous work 
by Y. Kawamata, H. Esnault and E. Viehweg. Applied to log 
varieties, our vanishing theorem is stronger
than Kawamata-Viehweg (or Nadel) vanishing.

We expect that normal quasi-log varieties are 
{\em equivalent} (cf. ~\ref{basicexmp}.1) to generalized 
log varieties, according to the Adjunction Conjecture.
We only have partial results in this direction 
(cf.~\ref{normal}, ~\ref{elt}, ~\ref{onadj}). One should 
also note that if the Adjunction Conjecture holds, the 
X-method works inductively in the category of log varieties, 
as long as we restrict to normal lc centers.

Finally, for technical reasons, we require that our varieties
with normal crossings singularities are globally embedded as 
hypersurfaces. This is enough for applications to generalized
log varieties, but we expect that this extra assumption is 
not necessary (see ~\ref{wsh}).

\begin{ack}
I would like to thank Professor Vyacheslav V. Shokurov for 
useful discussions. The motivation behind this work is his idea, 
that log varieties and their LCS loci should be treated on an equal
footing. I would also like to thank Professor Yujiro Kawamata 
for useful comments.

The author is currently a European Community Marie Curie Fellow.
The core of this work has been supported by the Japan Society for 
Promotion of Science, and partial support
by NSF Grant DMS-9800807 was received at an initial stage.

\end{ack}


\section{Preliminary}


A {\em variety} is a scheme of finite type, defined over an 
algebraically closed field $k$ of characteristic zero. We 
denote by $\Div(X)$ the abelian group of Cartier divisors 
of $X$. A {\em $K$-Cartier divisor} on $X$ is an element of 
$\Div(X)_K:=\Div(X)\otimes_\Z K$, for $K\in \{\Z,\Q,\R\}$.

Let $\pi:X\to S$ be a proper morphism of varieties.
We denote by $Z_1(X/S)$ the abelian group generated 
by proper integral curves in $X$ mapped to points by 
$\pi$. The natural pairing $\Div(X)\times Z_1(X/S)\to \Z$
induces, via numerical equivalence and tensoring 
with $\R$, a perfect pairing of finite dimensional 
$\R$-vector spaces
$
N^1(X/S)\times N_1(X/S) \to \R.
$
We denote by $NE(X/S)\subset N_1(X/S)$ the cone 
generated by proper integral curves in $X$ mapped to 
points by $\pi$, and by $\overline{NE}(X/S)$ its closure 
in the real topology. The dual of $NE(X/S)$ in $N^1(X/S)$
is called the {\em relatively nef cone}. 
The {\em relatively ample cone} $\Amp(X/S)$ is the cone of
$N^1(X/S)$ generated by classes of relatively ample Cartier 
divisors (if any).
A $K$-Cartier divisor $D$ is {\em relatively nef (ample)} if 
its class in $N^1(X/S)$ belongs to the relatively nef (ample) 
cone. 
If $X/S$ is projective, S. Kleiman proved that the relatively 
ample cone is the interior of the relatively nef cone. 
In particular, a $K$-Cartier divisor $D$ is relatively ample
if and only if  $(D \cdot z) > 0$ for all 
$z \in \overline{NE}(X/S)\setminus \{0\}$.

An $\R$-Cartier divisor $D$ is {\em relatively semi-ample} 
if $D\sim_\R f^*H$, where $f:X/S\to Y/S$ is a proper morphism, 
and $H$ is a relatively ample $\R$-Cartier divisor.
If $D\in \Div(X)_\Q$, this is equivalent to the surjectivity
of the natural map $\pi^*\pi_* \cO_X(mD) \to \cO_X(mD)$ for
some large and divisible positive integer $m$.

An open subset $U\subseteq X$ is called {\em big}
if $X\setminus U$ has codimension at least two in $X$.


\section{Normal crossings pairs}


\begin{defn}\label{multi}
A variety $X$ has {\em multicrossings} singularities if for every 
closed point $x \in X$, there exist integers $N, l$, subsets 
$I_1,\ldots,I_l$ of $\{0,\ldots,N\}$, and an isomorphism of 
complete local rings
$$
\cO^{\wedge}_{X,x} \isoto 
  \frac{k[[x_0,\ldots,x_N]]}{ (\prod_{i\in I_1} x_i,
  \ldots,\prod_{i\in I_l} x_i)}
$$
If $l=1$ for every $x\in X$, we say that $X$ has 
{\em normal crossings} singularities. Furthemore, if each 
irreducible component of $X$ is non-singular, we say that 
$X$ is a {\em simple multicrossings (normal crossings)} 
variety.
\end{defn}

For a scheme $X$, we denote by $\epsilon :X_\bullet \to X$ 
the associated simplicial scheme 
$((X_0/X)^{\Delta_n} \to X)_{n\ge 0}$. Here 
$\epsilon=\{\epsilon_n\}$, where $\epsilon_0:X_0\to X$ is 
the normalization, and $\epsilon_n$ is the natural projection. 
The simplicial maps are $\delta_i: X_{n+1}\to X_n, \ x_0 
\times \cdots \times x_{n+1} 
\mapsto x_0 \times \cdots \hat{x_i} \cdots \times x_n$ and
$s_i: X_n \to X_{n+1}, \ x_0 \times \cdots \times x_n 
\mapsto x_0 \times \cdots \times 
x_i \times x_i \times x_{i+1} \cdots \times x_n$.
This is a proper hypercovering \cite{Hodge3, plurisyst}.
A {\em strata} of $X$ is by definition the image on $X$ 
of some irreducible component of $X_\bullet$. 

\begin{lem}\label{basic1} The following hold for a variety 
$X$ with multicrossings singularities:
\begin{itemize}
     \item[(i)] The associated hypercovering 
	 $\epsilon:X_\bullet \to X$ 
	 is proper, smooth and of cohomological descent with respect 
	 to locally free sheaves on $X$. 
     \item[(ii)] We have an isomorphism of functors 
	 $Hom(X,\cdot)\isoto Hom(X_\bullet,\cdot)$.     
	 \item[(iii)] $X$ has seminormal singularities.
	 \item[(iv)] If $X$ is a simple multicrossings variety, 
	 each strata of $X$ is non-singular.
\end{itemize}
\end{lem}

\begin{proof} 
(i) Each $\epsilon_n$ is a finite map, thus proper. It is also 
easy to see that each $X_n$ is non-singular: in the notations 
of Definition~\ref{multi}, for $\alpha\in I_1\times \cdots 
\times I_l$ denote $\{\alpha\}=\{\alpha_1,\ldots,\alpha_l\} 
\subset \Delta_N$. Also, denote by $J$ the elements of 
$I_1\times \cdots \times I_l$ which are minimal with
respect to the partial order $\alpha\le \beta$ if and only if
$\{\alpha\} \subseteq \{\beta\}$. Then, at the complete 
local rings level, $X_n$ is the spectrum of 
$$
\sum_{q:\Delta_n \to J}
 \frac{k[[x_0,\ldots,x_N]]}{(x_i: i\in \{q(0)\} 
 \cup \cdots \cup  \{q(n)\} )}
$$
Finally, cohomological descent for a locally free sheaf $F$ 
on $X$ means that
$
F\isoto R^\bull \epsilon_*(\epsilon^*F).
$
Since $\epsilon$ is finite, it is enough to show that the 
natural map $\cO_X \to \epsilon_* \cO_{X_\bullet}$ is an 
isomorphism. This is a local statement, and it can be 
checked as in ~\cite[4.1]{plurisyst}.

(ii) A morphism $f:X\to Y$ induces $f:X_\bullet\to Y$ with 
components $f_n=f\circ \epsilon_n$. 
Conversely, let $f:X_\bullet\to Y$ be a morphism. The induced 
map $f:X\to Y$ is defined set-theoretically by
$
f(x):=f_0(\epsilon_0^{-1}(x)).
$
This map is well defined since any two points in the 
fiber of $\epsilon_0$ are the images of some point on 
some $X_n$ under different compositions of $\delta_i$'s.
Moreover, $f$ is a morphism since for every $h\in \cO_Y$, 
$f_0^*(h)\in \cO_{X_0}$ takes the same value on the 
glueing data, thus belongs to $\cO_X\subset \cO_{X_0}$.

(iii) See \cite{LCS}.

(iv) The normalization $\epsilon_0$ is a disjoint union 
of embeddings. Therefore the same holds for 
$\epsilon_n,\ n\ge 1$. Each $X_n$ is smooth since $X$ 
has multicrossings singularities, hence all strata are smooth. 
The strata are the components of the intersections of 
irreducible components of $X$, in this case.
\end{proof}

Let $X$ be a variety with multicrossings singularities. A 
Cartier divisor $D$ on $X$ is called {\em permissible} 
if it induces a Cartier divisor $D^\bullet$ on $X_\bullet$, i.e. 
$D^n=\epsilon_n^* D$ is a Cartier divisor on $X_n$, 
for every $n$ (equivalently, $D$ contains no strata of $X$
in its support).
We say that $D$ is a {\em multicrossings divisor} 
on $X$ if, in the notations of Definition~\ref{multi}, we have
$$
\cO^{\wedge}_{D,x} \isoto 
\frac{k[[x_0,\ldots,x_N]] }{(\prod_{i\in I_1} 
x_i,\ldots,\prod_{i\in I_l} x_i, 
\prod_{i\in I'} x_i)},
$$
where $I'\subset \Delta_N$ and $I'\cap \bigcup_{j=1}^l 
I_j=\emptyset$.
We denote by $\Div_0(X)$ the free abelian group generated
by all permissible Cartier divisors on $X$. A  permissible 
$K$-divisor on $X$ is an element of $\Div_0(X) \otimes_\Z K$, 
for $K\in \{\Z,\Q,\R\}$. 
For a permissible $K$-divisor $D=\sum_i d_i D_i$, its 
{\em reduced part} is $\sum_{d_i=1}D_i$. We denote 
$D^{>1}=\sum_{d_i>1}d_i D_i$ and $D^{<1}=\sum_{d_i<1}d_i D_i$. 
We say that $D$ is a {\em boundary} ({\em sub-boundary}) if 
$0\le d_i\le 1 \ \forall i$ ($d_i\le 1 \ \forall i$).

\begin{defn}
A {\em multicrossings pair} $(X,B)$ is a multicrossings variety $X$ 
endowed with a permissible $\R$-divisor $B$, whose support 
is a multicrossings divisor on $X$. If $X$ has normal crossings 
singularities, we say that $(X,B)$ is a {\em normal crossings pair}.

A {\em strata of} $(X,B)$ is a strata of either $X$, or the 
reduced part of $B$. Equivalently, the strata are the images 
of strata of the log-nonsingular pairs $\{(X_n,B^n)\}_{n\ge 0}$. 
For instance, the maximal strata of $(X,B)$ are the irreducible 
components of $X$.
\end{defn}

\begin{rem} Compared with the {\em generalized normal crossings 
varieties} introduced by Y. Kawamata ~\cite{plurisyst}, the 
ambient space $X$ of a normal crossings pair has generalized 
normal crossings singularities, but $B$ has arbitrary 
coefficients in our case. 
\end{rem}

\begin{lem} The following properties hold for a multicrossings 
pair $(X,B)$:
\begin{itemize}
\item[(i)] Each strata is irreducible, with multicrossings 
  singularities. A strata which is minimal (with respect to inclusion) 
  is non-singular.
\item[(ii)] There are only finitely many strata.
\item[(iii)] The non-empty intersection of any two strata 
  is a union of strata. In particular, minimal strata are 
  mutually disjoint.
\end{itemize}
\end{lem}
We say that a permissible divisor $D$ has 
{\em multicrossings support on} $(X,B)$ if it contains
no strata of $(X,B)$ and both $D$ and its restriction to the
reduced part of $B$ have multicrossings support. A variety 
with normal crossings $X$ is locally complete intersection, 
so it has an invertible dualizing sheaf $\cO_X(K)$. The 
{\em canonical divisor} $K\in \Div(X)$ is well defined up to 
linear equivalence.

\begin{rem} [D\'evissage] \label{adj} 
Let $(X,B)$ be a normal crossings pair, and let $Y$ be a union 
of irreducible components of $X$. Denote by $X'$ the union of 
the other irreducible components of $X$, and write
$B_Y=B|_Y+{X'}|_Y$, $B_{X'}=Y|_{X'}+B|_{X'}$. 
Then the following hold:
\begin{itemize}
 \item[(i)] $(Y,B_Y)$ and $(X',B_{X'})$ are normal crossings pairs.    
 \item[(ii)] $(K+B)|_Y=K_Y+B_Y$ and $(K+B)|_{X'}=K_{X'}+B_{X'}$.
 \item[(iii)] $\cI_{Y,X} \simeq j_* \cO_{X'}(-Y|_{X'})$, where 
 $j:X'\to X$ is the inclusion.
 \end{itemize}
 In particular, let $L$ be a Cartier divisor on $X$ such that 
 $L= K+B+H$. Denote $L'=L|_{X'}-Y|_{X'}$, 
 so that $L'= K_{X'}+B|_{X'}+H|_{X'}$. Then we have a short 
 exact sequence 
 $0 \to j_*\cO_{X'}(L') \to \cO_X(L) \to \cO_Y(L|_Y) \to 0$. 
\end{rem}

\begin{defn} We say that a normal crossings pair $(X,B)$ is
{\em embedded} if there exists a closed embedding $j:X \to M$, 
where $M$ is a non-singular variety of dimension $\dim(X)+1$.
\end{defn}

Let $(X,B)$ be an embedded normal crossings pair, and let
$C$ be a nonsingular strata. The {\em embedded log 
transformation of $(X,B)$ in $C$}, denoted 
$\sigma:(Y,B_Y)\to (X,B)$, is defined as follows: let 
$X \subset M$ be an embedding of $X$ as a hypersurface in a 
nonsingular ambient space $M$.
We denote by $Y$ the reduced structure of the total 
transform of $X$ in the blow-up of $M$ in $C$. 
The morphism $\sigma:Y\to X$ is projective, $Y$ has 
normal crossings singularities, the formula 
$\sigma^*(K+B)=K_Y+B_Y$ defines a divisor $B_Y$ on $Y$, 
and the following properties hold:
\begin{itemize}
\item[(i)] $(Y,B_Y)$ is an embedded normal crossings pair.
\item[(ii)] The strata of $(X,B)$ are exactly the images of 
    the strata of $(Y,B_Y)$.
\item[(iii)] $\cO_X \isoto R^\bull \sigma_*\cO_Y$.
\item[(iv)] $\sigma^{-1}(C)$ is a maximal strata of $(Y,B_Y)$.
\end{itemize}

\begin{prop}\label{multi_lcs} Let $X' \subset X$ 
be the union of some strata of an embedded normal crossings 
pair $(X,B)$. Then there exists an embedded normal crossings 
pair $(Y,B_Y)$, and a projective morphism $f:Y\to X$ such that:
 \begin{itemize}
   \item[(i)]  $\cO_X\isoto R^\bull f_*\cO_Y$.
   \item[(ii)] $f^*(K+B)=K_Y+B_Y$.
   \item[(iii)] The strata of $(X,B)$ are exactly the images 
       of the strata of $(Y,B_Y)$.
   \item[(iv)] $f^{-1}(X')$ is a union of maximal strata of 
             $(Y,B_Y)$.
 \end{itemize}
\end{prop}

\begin{proof} First, we may assume that each strata of 
$(X,B)$ is nonsingular. Indeed, after a finite number of embedded
log transformations of $X$ in its minimal strata, each irreducible 
component of $X$ is nonsingular in the minimal strata of $X$, 
i.e. $X$ has simple normal crossings. Similarly, the reduced 
part of $B$ becomes simple multicrossings after a finite 
sequence of embedded log transformations of $(X,B)$ in minimal 
strata of $B$.

Once each strata of $(X,B)$ is nonsingular, we reach the 
conclusion after a finite number of embedded log transformations 
of $(X,B)$ in the irreducible components of $X'$.
\end{proof}

\begin{rem} \label{wsh}
The embedded hypothesis is used to prove ~\ref{multi_lcs}, 
and to resolve singularities of permissible subvarieties of 
a variety with normal crossings. Once the latter has been 
established, we expect our results to work for abstract normal 
crossings pairs.
\end{rem}


\section{Vanishing theorems}


We extend the vanishing and torsion freeness theorems of 
J. Koll\'ar ~\cite{higherimag} to normal crossings pairs. 
The proof is based on logarithmic De Rham complexes, and 
we follow closely the presentation of \cite{lectures}. 
See also ~\cite{plurisyst}.

\begin{thm} \label{tako} Assume $(X,B)$ is an embedded normal 
crossings pair such that $X$ is a proper variety and $B$ 
is a boundary. Let $L$ be a Cartier divisor on $X$ and let 
$D$ be an effective Cartier divisor, permissible with respect to $(X,B)$, 
with the following properties:
\begin{itemize}
   \item[(i)] $L\sim_\R K+B+H$.
   \item[(ii)] $H \in \Div(X)_\R$ is semi-ample.
   \item[(iii)] $tH \sim_\R D+D'$ for some positive real
       number $t$, and for some effective $\R$-Cartier 
       divisor $D'$, permissible with respect to to $(X,B)$.
\end{itemize}
Then the natural maps $H^q(X,\cO_X(L))\to H^q(X,\cO_X(L+D))$ 
are injective for all $q$.
\end{thm}

\begin{proof} Blowing up $X$ and incorporating the negative
part of $B$ into the pullback of $L$, we may assume that both 
$(X,B)$ and $D+D'$ have normal crossings support. Furthemore,
we may assume $H=a D+a' D'$, where $a>0, a'\ge 0$, 
and $B'=B+a D+a' D'$ is a boundary with 
$\lfloor B' \rfloor =\lfloor B\rfloor$.

We have $L\sim_\R K+B'$. Since $L,K$ are integral divisors,
the set of boundaries having the same support and reduced part 
as $B'$ and satisfying the above equality,
form a rational polyhedra. After a perturbation of its fractional
part, we may assume that $B'$ is rational. In particular, 
$T=-L+K+B'$ is a $\Q$-Cartier divisor and $\nu T\sim 0$ for
some positive integer $\nu$. Assume that $\nu$ is minimal
with this property. Denote $\cE=\cO_X(-L+K)$, and let $R$ be the 
support of $B'$.

Let $X_\bullet \to X$ be the associated smooth, proper 
hypercovering. By Serre duality and cohomological descent, 
we have to check the surjectivity of the maps 
$$
H^q(X_\bullet,\cE^\bullet (-D^\bullet)) \to 
H^q(X_\bullet,\cE^\bullet).
$$
We use the following commutative diagram:
\[ \xymatrix{
H^q(X_\bullet,\cE^\bullet (-D^\bullet))  \ar[r] & 
H^q(X_\bullet,\cE^\bullet)    \\
{\bf H}^q(X_\bullet,\Omega_{X_\bullet}^\bull (\mbox{log } 
R^\bullet) 
\otimes \cE^\bullet(-D^\bullet))  \ar[u] \ar[r]^{\alpha} & 
{\bf H}^q(X_\bullet,\Omega_{X_\bullet}^\bull (\mbox{log } 
R^\bullet) \otimes \cE^\bullet)  \ar[u]_\beta
} \]
Since $-L+K =\lfloor T \rfloor -\lfloor B \rfloor$, the 
restriction of $\cE^\bullet$ to each component of $X_\bullet$
admits a logarithmic connection with poles along $R^\bullet$, 
whose residues along the components of $D^\bullet$ belong to the 
interval $(0,1)$ \cite[3.2]{lectures}. 
By \cite[4.3]{lectures}, the map 
$$
\Omega_{X_\bullet}^\bull (\mbox{log } R^\bullet) \otimes 
        \cE^\bullet(-D^\bullet) \to
\Omega_{X_\bullet}^\bull (\mbox{log } R^\bullet) 
	 \otimes \cE^\bullet
$$
is a quasi-isomorphism componentwise, thus it is a 
quasi-isomorphism of simplicial complexes. Therefore 
$\alpha$ is an isomorphism. 

Let $\pi:Y_\bullet \to X_\bullet$ be the cyclic cover of 
degree $\nu$ corresponding to the torsion divisor $T^\bullet$. 
By \cite{Hodge3}, the spectral sequence
$$
E^{pq}_1=H^q(Y_\bullet,\Omega_{Y_\bullet}^p (\mbox{log } 
R^\bullet)) \Longrightarrow
{\bf H}^{p+q}(Y_\bullet, \Omega_{Y_\bullet}^\bull 
(\mbox{log } R^\bullet))
$$
degenerates. Since $\cE^\bullet$ is a direct summand of 
$\pi_*\Omega^\bull_{Y_\bullet}(\mbox{log } R^\bullet)$,
the spectral sequence
$$
E^{pq}_1=H^q(X_\bullet,\Omega_{X_\bullet}^p (\mbox{log } R^\bullet) 
           \otimes \cE^\bullet)
\Longrightarrow
{\bf H}^{p+q}(X_\bullet, \Omega_{X_\bullet}^\bull (\mbox{log } 
R^\bullet) \otimes \cE^\bullet)
$$
degenerates as well. Therefore $\beta$ is surjective.
\end{proof}

\begin{thm}\label{tf} 
Let $(Y,B)$ be an embedded normal crossings pair, and 
assume that $B$ is a boundary. Let $f:Y\to X$ be a proper 
morphism, and let $L$ be a Cartier divisor on $Y$ such that 
$H=L-(K+B)$ is $f$-semi-ample. Then:
\begin{itemize}
 \item[(i)] Every non-zero local section of $R^qf_*\cO_Y(L)$ 
  contains in its support the $f$-image of some strata 
  of $(Y,B)$.
 \item[(ii)] Let $\pi:X\to S$ be a projective morphism 
   and assume $H\sim_\R f^*H'$ for some $\pi$-ample 
   $\R$-Cartier divisor $H'$ on $X$. 
   Then $R^qf_*\cO_Y(L)$ is $\pi_*$-acyclic.
  \end{itemize}
\end{thm}

\begin{proof}
(i) The conclusion is local, so we may shrink $X$ to 
an affine open subset and compactify it afterwards, 
so that $X$ is projective, $Y$ is proper and $H$ is 
semi-ample. If $R^qf_*\cO_Y(L)$ admits a local section 
whose support does not contain any image of the $(Y,B)$ 
strata, one can find a very ample divisor $A$ such that:
 \begin{itemize}
 \item[-] 
   $H^0(X,R^qf_*\cO_Y(L))\to H^0(X,R^qf_*\cO_Y(L)
   \otimes \cO_X(A))$ 
    is not injective.
 \item[-] $f^*A$ is a permissible multicrossings 
   divisor on $(Y,B)$.
 \item[-] The Leray spectral sequence of $L+f^*A$ with 
   respect to $f$ degenerates. 
  \end{itemize}
Replacing $L$ by $L+f^*A$ if necessary, we may also assume that 
$H-f^*A$ is semiample. The degeneration of the Leray spectral 
sequence implies that the map 
$H^q(Y,\cO_Y(L)) \to H^q(Y,\cO_Y(L+f^*A))$ is not 
injective, which contradicts Theorem~\ref{tako}.

(ii) Assume $\dim S=0$, and let $H=f^*H_X$. 
If $X$ has positive dimension, one can find a divisor $A$
in some large, divisible multiple of $H$, such that its pullback 
$A'=f^*A$ is a permissible multicrossings divisor on $(Y,B)$, and 
$R^qf_*\cO_Y(L+A')$ is $\pi_*$-acyclic for all $q$. By (i), 
we have short exact sequences:
$$
0 \to R^qf_*\cO_Y(L) \to 
R^qf_*\cO_Y(L+A')\to R^qf_*\cO_{A'}(L+A') \to 0
$$
$R^qf_*\cO_Y(L+A')$ is $\pi_*$-acyclic by assumption, while 
$R^qf_*\cO_{A'}(L+A')$ is $\pi_*$-acyclic by induction on 
$X$. Therefore $E^{p,q}_2=0$ for $p\ge 2$ in the 
following commutative diagram of spectral sequences:
\[ \xymatrix{
    E^{p,q}_2= R^p\pi_* R^qf_*\cO_Y(L)  
    \ar[d]^{\varphi^{p,q}} &\Longrightarrow &  
    R^{p+q}(\pi \circ f)_*\cO_Y(L) 
    \ar[d]^{\varphi^{p+q}}   \\
    \bar{E}^{p,q}_2= R^p\pi_* R^qf_*\cO_Y(L+A')  
    & \Longrightarrow  & R^{p+q}(\pi \circ f)_*\cO_Y(L+A')
} \]
Since $E^{1,q}_2 \to R^{1+q}(\pi \circ f)_*\cO_Y(L)$ is 
injective, $\varphi^{1+q}$ is injective by Theorem~\ref{tako}, 
and $\bar{E}^{1,q}=0$ by assumption, we obtain $E^{1,q}_2=0$. 

Assume now the $S$ is affine of positive dimension, and 
$\pi\circ f$ surjects $Y$ onto $S$. We use induction on the 
dimension of $S$. 
\begin{itemize}
   \item[a)] Assume that each strata of $(Y,B)$ dominates a 
generic point of $S$. From the case $\dim S=0$, 
$R^p\pi_*R^qf_*\cO_Y(L) (p>0)$ does not contain any generic 
point of $S$ in its support. Therefore there exists a 
general hyperplane section $A$ of $S$, containing the 
support of all these sheaves, such that its pullback 
$A'$ on $Y$ is a multicrossings divisor on $(Y,B)$.
The argument in (i) shows that $R^qf_*\cO_Y(L)$ is 
$\pi_*$-acyclic, except that $\varphi^{p+q}$ is injective 
by (i) now, and $R^p\pi_*R^qf_*\cO_Y(L) \otimes \cO_S(A)$ is 
zero by the choice of $A$.
   \item[b)] Let $Y'$ be the union of all strata of 
$(Y,B)$ which is not mapped onto generic points of $S$. 
After a sequence of embedded log transformations, we 
may assume that $Y'$ is a union of irreducible components 
of $Y$. By $(i)$, we have exact sequences
$$
0 \to R^q f_*(\cI_{Y'}(L)) \to
R^qf_*\cO_Y(L) \to R^qf_*\cO_{Y'}(L) \to 0.
$$
>From Remark~\ref{adj}, 
$R^q f_*(\cI_{Y'}(L)) \isoto R^q f_*\cO_{Y''}(L'')$, 
where $L''=K_{Y''}+B|_{Y''}+f^*H$. 
The pair $(Y'', B|_{Y''})$ satisfies the hypothesis 
in $a)$, hence the first term is $\pi_*$-acyclic. The third 
is $\pi_*$-acyclic by induction, thus $R^qf_*\cO_Y(L)$ 
is $\pi_*$-acyclic.	   
    \end{itemize}
\end{proof}


\section{Quasi-log varieties}


\begin{defn} \label{qlv} 
 A {\em quasi-log variety} is a scheme $X$ 
 endowed with an $\R$-Cartier divisor $\omega$, a proper 
 closed subscheme $X_{-\infty} \subset X$, and a finite
 collection $\{C\}$ of reduced and irreducible subvarieties
 of $X$ such that there exists a proper morphism 
 $f:(Y,B_Y)\to X$ from an embedded normal crossings pair 
 satisfying the following properties:
 \begin{itemize}
 \item[(1)] $f^*\omega \sim_\R K_Y+B_Y$.
  \item[(2)] The natural map 
  $\cO_X\to f_*\cO_Y(\lceil -(B_Y^{<1}) \rceil)$ 
  induces an isomorphism 
  $\cI_{X_{-\infty}} \to f_*\cO_Y(\lceil -(B_Y^{<1}) 
  \rceil - \lfloor B_Y^{>1} \rfloor)$. 	                  
  \item[(3)] The collection of subvarieties $\{C\}$ 
   coincides with the images of $(X,B)$-strata which are
   not included in $X_{-\infty}$.
  \end{itemize}
 We use the following terminology: the subvarieties $C$ 
 are the {\em qlc centres} of $X$, $X_{-\infty}$ is  
 the {\em non-qlog canonical locus} of $X$, and 
 $f:(Y,B)\to X$ is a {\em quasi-log resolution} of $X$.
 We say that $X$ has {\em qlog canonical singularities} if 
 $X_{-\infty}=\emptyset$. Note that a quasi-log variety $X$ 
 is the union of its qlc centers and $X_{-\infty}$.
 A {\em relative quasi-log variety} $X/S$ is a quasi-log 
 variety $X$ endowed with a proper morphism $\pi:X\to S$.

 For simplicity, we will refer to a quasi-log variety
 as $X$ or $(X,\omega)$. 
 \end{defn}
 
 \begin{rem}
 \begin{itemize}

\item[(i)] $X$ has qlog canonical singularities if 
and only if $B$ is a sub-boundary. Indeed, the commutative 
diagram
\[ \xymatrix{
0 \ar[r] & f_*\cI_N \ar[r]  & f_*\cO_Y \ar[r] & f_*\cO_N & \\
0 \ar[r] & \cI_{X_{-\infty}} \ar[r] \ar[u]_{\simeq} & \cO_X 
\ar[r] \ar[u]  &  \cO_{X_{-\infty}} \ar[r]\ar[u]  &  0
 } \]
implies that $X_{-\infty} \cap f(Y)=f(N)$, where 
$N=\lfloor B_Y^{>1} \rfloor$. Note that 
$X_{-\infty}\subsetneq X$ by assumption, but $X_{-\infty}$ 
may contain irreducible components of $X$. Also, $f$ may not 
be surjective (cf. ~\ref{basicexmp}.4).
 
\item[(ii)] If $B$ is a sub-boundary, ~\ref{qlv}.2
says that the natural morphism 
$\cO_X \to f_*\cO_Y(\lceil -(B_Y^{<1}) \rceil)$ 
is an isomorphism.  In particular, $f$ is a 
surjective map with connected fibers. Furthemore,
$X$ is seminormal by \cite{LCS}. 
In general, the same holds over the open subset
of qlog-canonical singularities $U=X\setminus 
X_{-\infty}$. 

\item[(iii)] The {\em quasi-log canonical class} $\omega$
  is defined up to $\R$-linear equivalence. This is more
  general than the case of generalized log varieties, where 
  the log canonical class $K+B$ is defined up to linear 
  equivalence.

\item[(iv)] The quasi-log resolution plays a role similar to 
  a log resolution. Embedded log transformations of $(Y,B_Y)$,
  or blow-ups of $Y$ in centers which contain no 
  $(Y,B_Y)$-strata, leave the quasi-log structure on $X$ 
  invariant. 
  Furthermore, we may sligthly perturb the non-reduced 
  components of $B$.
  In particular, if $\omega$ is a $\Q$-divisor, we may 
  assume that $B$ is a $\Q$-divisor.

  \begin{proof} We check the invariance of the structure
  under permissible blow-ups (for embedded log transformations is 
  easier). The blow-ups do not introduce new $(Y,B)$-strata
  so we only need to check the invariance of the ideal sheaf
  in ~\ref{qlv}.2. By cohomological descend, we may assume 
  $Y$ is non-singular and $B_Y$ is a divisor with normal 
  crossings support. Assume $\sigma:(Y',B_{Y'}) \to (Y,B_Y)$ 
  is a crepant log non-singular model. Denote 
  $\Delta=B_Y-\lfloor B_Y\rfloor$, let $R$ be the reduced part 
  of $B_Y$, and define $\Delta'$ and $R'$ similarly. Note the 
  identity 
  $$
  \lceil -(B_Y^{<1}) \rceil - \lfloor B_Y^{>1} \rfloor
  =\lceil -B_Y \rceil +R.
  $$
  We have
  $
  (\lceil -B_{Y'} \rceil +R)-\sigma^*(\lceil -B_{Y'} \rceil +R)
  =K_{Y'}+\Delta'+R'-f^*(K_Y+\Delta+R).
  $
  It is enough to show that the right hand side is effective.
  Assume that is negative in some divisor $E$. Its coefficient
  $\mult_E(\Delta'+R')+a(E;\Delta+R)-1$ is integral, hence
  $\mult_E(\Delta'+R')+a(E;\Delta+R)\le 0$ (here $a(E;\Delta+R)$
  is the log discrepancy of $E$ with respect to $(Y,\Delta+R)$).
  Therefore $\mult_E(\Delta'+R')=0$ and $a(E;\Delta+R)=0$. 
  The latter implies that $c_Y(E)$ is a strata of $R$, hence
  we also have $a(E;B_Y)=0$ by the normal crossings assumption.
  Equivalently, $\mult_E(R')=1$. Contradiction. 
  \end{proof}
 \end{itemize}
 \end{rem}  
 
\begin{exmp} \label{basicexmp}
 \begin{itemize}
  \item[1.] Any generalized log variety $(X,B)$ is a 
  a quasi-log variety: let $\omega$ be any
  $\R$-Cartier divisor such that $\omega \sim_\R K+B$, 
  and let $X_{-\infty}$ be the locus where $(X,B)$ does not
  have log canonical singularities (with the induced closed
  subscheme structure). A quasi-log resolution is a log 
  resolution. The qlc centers are exactly the subvarieties
  $C$ of $X$ such that $(X,B)$ has zero log discrepany
  in the generic point of $C$. With the exception of $X$
  (which is a qlc centre), the qlc centres of $(X, \omega)$ 
  are exactly the {\em lc centres} of Y. Kawamata ~\cite{fujita34} 
  which are not included in $(X,B)_{-\infty}$. This is natural, 
  since we do not expect any adjunction on lc centres along which
  $(X,B)$ does not have log canonical singularities.

  Conversely, if $Y$ is non-singular, $f$ is birational 
  and $X$ is normal, then $X$ is associated (equivalent)
  to a generalized log variety as above. Indeed, the 
  corresponding generalized log variety is $(X,f_*B_Y)$.

\item[2.] Let $(Y,B_Y)$ be a proper log variety
  such that $K_Y+B_Y$ is nef. The 
  Abundance Conjecture predicts the existence of a proper 
  morphism $f:Y\to X$ to a projective variety $X$ such that 
  $K_Y+B_Y\sim_\R f^*H$ for some ample divisor
  $H\in \Div(X)_\R$. Then $X$ is a quasi-log
  variety with qlog canonical singularities, with 
  $\omega \sim_\R H$ and quasi-log resolution $f$.

\item[3.] Let $(\bar{X},\bar{B})$ be a generalized log variety, 
and assume that $X=LCS(\bar{X},\bar{B})$ intersects the 
open subset on which $(\bar{X},\bar{B})$ has log 
canonical singularities. Then $X$ is a quasi-log variety,
where $\omega \sim_\R (K_{\bar{X}}+\bar{B})|_X$ and 
$X_{-\infty}=(\bar{X},\bar{B})_{-\infty}$.
A quasi-log resolution of $X$ is induced by restricting 
to the reduced part of the boundary on a log resolution of 
$(\bar{X},\bar{B})$:
\[ \xymatrix{
(Y,B_Y) \ar[d]_f \ar[r]  & (\bar{Y},\bar{B}) \ar[d]_\mu\\
X    \ar[r]  & (\bar{X},\bar{B}) }\]
Here $K_{\bar{Y}}+
\bar{B}=\mu^*(K_{\bar{X}}+\bar{B})$, $Y$ is the reduced 
part of $\bar{B}$, and $B_Y=(\bar{B}-Y)|_Y$.

\item[4.] Let $X$ be a divisor with normal crossings 
support in a non-singular variety $\bar{X}$, and assume 
that $Y$, the reduced part of $X$, is non-empty. 
Then $X$ is a quasi-log variety, where 
$\omega \sim_\R (K_{\bar{X}}+X)|_X$, and $X_{-\infty}$ is the union 
of non-reduced components of $X$. A quasi-log resolution 
is $f:(Y,B_Y)\to X$, where $B_Y$ is defined by the 
adjunction formula $K_Y+B_Y=(K_{\bar{X}}+X)|_Y$.
   \end{itemize}
 \end{exmp}
 
 \begin{thm}[Adjunction \& Vanishing] \label{adj_van}
 Let $X$ be a quasi-log variety, and let $X'$ be 
 the union of $X_{-\infty}$ with a (possibly empty) 
 union of some qlc centers of $X$. 
 \begin{itemize}
 \item[(i)] 
   Assume $X'\ne X_{-\infty}$. Then $X'$ is a quasi-log 
  variety, with $\omega'=\omega|_{X'}$ and 
  $X'_{-\infty}=X_{-\infty}$. 
  Moreover, the qlc centers of $X'$ are exactly the qlc 
  centers of $X$ which are included in $X'$.
 \item[(ii)] 
   Assume $X/S$ is projective and let $L \in \Div(X)$ 
   such that $L-\omega$ is $\pi$-ample. Then
  $\cI_{X'}\otimes \cO_X(L)$ is $\pi_*$-acyclic.
     \end{itemize}
 \end{thm}
 
 \begin{proof} (i) After embedded log transformations, 
 we may assume that the union of all strata of $(Y,B_Y)$ mapped 
into $X'$, which we denote $Y'$, is a union of irreducible 
components of $Y$.
 Define $B_{Y'}$ by $(K_Y+B_Y)|_{Y'}=K_{Y'}+B_{Y'}$. We claim that 
 $f:(Y',B_{Y'})\to X'$ is a quasi-log resolution. The adjunction 
formula is clear, so we just check the second property.
 Denote $A=\lceil -(B_Y^{<1}) \rceil$ and $N= \lfloor B_Y^{>1} \rfloor$.
 Let $Y''$ be the subscheme of $Y$ whose ideal sheaf $\cI$ is defined 
 by the exact sequence
 $$
 0 \to \cI \to \cO_Y(-N) \to  \cO_{Y'}(-N) \to 0
 $$
 The ideal of the subscheme $X'$ is the unique ideal sheaf 
 $\cI_{X'} \subset \cI_{X_{-\infty}}$ for which the induced map 
 $\cI_{X'} \to f_*\cI(A)$ is an isomorphism. Consider the following 
 commutative diagram:
      \[ \xymatrix{
0 \ar[r] & f_*\cI(A) \ar[r]\ar[d]^{=} & f_*\cO_Y(A-N) 
  \ar[r]\ar[d] & f_*\cO_{Y'}(A-N) \ar[r]\ar[d] & 0 \\
0 \ar[r] & f_*\cI(A) \ar[r]  & f_*\cO_Y(A) \ar[r]  & 
 f_*\cO_{Y''}(A)               &    \\
0 \ar[r] & \cI_{X'} \ar[r] \ar[u]_{\simeq} & \cO_X 
 \ar[r] \ar[u]  & \cO_{X'} \ar[r]\ar[u]  &  0
        } \]
 The map $f_*\cO_{Y'}(A-N)  \to f_*\cO_{Y''}(A)$ is injective by 
 the definition of $\cI$. Moreover, 
 $\cI(A) \simeq \cI_{Y'} \otimes \cO_Y(A-N)$ and $K_Y+B_Y \sim_\R 0/X$. 
 From the choice of $Y'$, we deduce from Theorem~\ref{tf}(i) that any 
 local section of $R^1f_*\cI(A)$ which is supported by $f(Y')$ is zero. 
 Therefore the top row is exact.  It is easy to see that
 $\cI_{X'_{-\infty}}:=\cI_{X_{-\infty}}/ \cI \to f_*\cO_{Y'}(A-N)$ 
 is an isomorphism. Finally, the characterization of the 
 qlc centers of $X'$ follows from the choice of $Y'$, 
 and the corresponding statement for $(Y',B_{Y'})$ and $(Y,B_Y)$.
 
 (ii) As in the proof of Theorem~\ref{tf}(ii.b), it follows from 
    Theorem~\ref{tf}(ii) that $f_*\cI(A) \otimes \cO_X(L)$ 
    is $\pi_*$-acyclic. 
 \end{proof}
 
 \begin{rem}\label{lift} 
 The above proof gives a commutative diagram of short exact 
 sequences:
 \[ \xymatrix{
  0 \ar[r]  &  \cI_{X'} \ar[r] \ar[d]_{=} & \cI_{X_{-\infty}} 
 \ar[r] \ar[d] &  \cI_{{X'}_{-\infty}} \ar[r] \ar[d]  &   0 \\
  0 \ar[r]  &  \cI_{X'} \ar[r]    & \cO_X   \ar[r]    &  
  \cO_{X'} \ar[r] &  0
  } \]
 Therefore we can lift global sections of 
 $\cO_X(L)$ or $\cI_{X_{-\infty}} \otimes \cO_X(L)$ 
 from $X'/S$ to $X/S$.
 \end{rem}
 
\begin{defn} The {\em LCS locus} of a quasi-log 
variety $X$, denoted $LCS(X)$, is $X_{-\infty}$ union 
with all qlc centers of $X$ which are not maximal with 
respect to the inclusion. The subscheme structure is 
defined as above, and we have a natural embedding
$X_{-\infty}\subseteq LCS(X)$.
\end{defn}

\begin{prop}\label{normal} 
Let $X$ be a quasi-log variety whose LCS locus 
is empty. Then $X$ is normal.
\end{prop}

\begin{proof} We may assume that $X$ is connected. 
Let $f:(Y,B_Y)\to X$ be a quasi-log resolution of $X$. 
By assumption, $B_Y$ is a sub-boundary, $f$ is surjective 
with connected fibers and each strata of $(Y,B_Y)$ 
dominates some irreducible component of $X$.
We first show that $X$ is irreducible. Indeed, let
$\{X_i\}$ be the irreducible components of $X$, and 
let $Y_i$ be the union of strata of $Y$ which 
dominate $X_i$. A non-empty intersection of two 
strata mapped on different components cannot 
dominate some component of $X$, thus $Y$ is the 
disjoint union of the closed subsets $Y_i$. But $Y$ 
is connected since $f$ has connected fibers, thus 
$X$ is irreducible.

Let $f_n:Y_n \to X$ be the induced morphisms. 
Then $Y_n=\sqcup_j Y_n^j$ is the disjoint union of 
its irreducible components, and $f_n=\sqcup f_n^j$. 
Each $f_n^j:Y_n^j \to X$ is dominant, thus factors 
through the normalization: $f_n^j=\nu\circ g_n^j$. 
The maps $\{g_n=\sqcup g_n^j\}_n$ glue to a 
morphism $g:Y_\bullet \to X^\nu$ which factors 
$f:Y_\bullet\to X$. This map extends to $Y$ according to 
Lemma~\ref{basic1}(ii).

Therefore $f$ factors through the normalization 
of $X$. Since $f$ has connected fibers and $X$ is 
seminormal, the normalization is an isomorphism.
\end{proof}

The following properties of qlc centers generalize 
~\cite[1.5, 1.6]{fujita34} (in particular, minimal
lc centers of log varieties have normal 
singularities):

\begin{prop}\label{conn} 
Assume $X$ is a quasi-log variety with qlog 
canonical singularities. The following hold:
\begin{itemize}
  \item[(i)] The intersection of two qlc centers 
    is a union of qlc centers.
  \item[(ii)] For any point $P\in X$, the set of all qlc 
    centers passing through $P$ has a unique minimal 
	element $W$. Moreover, $W$ is normal at $P$.
\end{itemize}
\end{prop}

\begin{proof}
\begin{itemize}

\item[(i)] Let $C_1,C_2$ be two qlc centers of $X$. 
Fixing $P\in C_1\cap C_2$, is enough to find an 
qlc center $C$ such that $P\in C\subset C_1\cap 
C_2$. $X'=C_1\cup C_2$ is a quasi-log variety 
with two irreducible components, hence it is not 
normal at $P$. By Proposition~\ref{normal},
$P\in LCS(X')$. Therefore there exists a qlc center 
$C\subset C_1$ with $\dim C<\dim C_1$ such that
$P\in C\cap C_2$. If $C\subset C_2$, we are done. 
Otherwise we repeat the argument with $C_1:=C$, 
and reach the conclusion in a finite number of steps.   
  
\item[(ii)] The uniqueness follows from (i), and the 
normality from Proposition~\ref{normal}.

\end{itemize}
\end{proof}

\begin{thm}(cf. ~\cite{adjkaw}) \label{elt} 
Let $(X/S,B)$ be a relative generalized log variety. 
Let $\nu:W \to X$ be the normalization of an irreducible 
component of $LCS(X,B)$, and assume $\nu(W)$ is an 
exceptional lc centre. The following hold:
\begin{itemize}
 \item[(i)] There exists a quasi-log structure on 
    $W$ such that $\omega \sim_\R \nu^*(K+B)$ and
    $LCS(W,\omega) \subseteq 
    \nu^{-1}((X,B)_{-\infty} \cup 
    \bigcup \{C \mbox{ lc centre}\ne \nu(W)\})$.
 \item[(ii)] Assume $H$ be a nef and big $\R$-divisor 
    on $W/S$. Then there exists a generalized log variety 
    structure $(W,B_W)$ on $W$ such that $\omega+H \sim_\R 
    K_W+B_W$ and $LCS(W,B_W) \subseteq LCS(W,\omega)$. 
\end{itemize}
\end{thm}

\begin{rem}\label{onadj} 
\begin{itemize}
\item[1.] This is a weak form of adjunction. We expect that the 
inclusion in (i) is an equality (we prove this on a big 
open subset of $W$). Furthermore, (ii) should hold in 
a stronger form: the quasi-log structure of $(W,\omega)$ 
is equivalent to the log structure of $(W,B_W)$.

\item[2.] $(X,B)$ induces a natural $\R$-b-divisor 
${\mathcal B}_{div}$ of $W$, called the {\em divisorial part 
of adjunction} (cf. ~\cite[\S 3]{thesis}), and the following 
inequality is expected to hold:
$$
{\mathcal A}(W,B_W) \le -{\mathcal B}_{div}
$$
If $\dim(X)\le 4$, this follows from ~\cite{PSh2}: there 
exists a birational model $W'/W$ such that 
$-{\mathcal B}_{div}={\mathcal A}(W',({\mathcal B}_{div})_{W'})$.
This implies the desired inequality if we choose a high enough 
model $W'/W$ in Step (ii) of the proof. 
\end{itemize}
\end{rem}

\begin{proof} 
\begin{itemize}
 \item[(i)] The lc centre being exceptional means
that among the valuations centered at $\nu(W)$ on $X$,
there exists a unique valuation $E$ with zero log
discrepancy with respect to $(X,B)$. 
Let $\mu:(Y,B_Y) \to (X,B)$ be a crepant log resolution 
such that $E$ is a divisor on $Y$.
We can write $B_Y=E+B'$, and set $B_E=B'|_E$ and 
$\omega=\nu^*(K+B)$. Since $f:E \to \nu(W)$ has connected 
general fibre, its Stein factorization is $g:E \to W$:
\[
\xymatrix{ (Y,B_Y) \ar[d]_\mu & (E,B_E) \ar[d]^g \ar[l]  \\
  (X,B) &  (W,\omega) \ar[l]^\nu
}
\]
We claim that $g$ defines a quasi-log structure on $W$.
Indeed, the crepant hypothesis is satisfied since 
$g^*\omega \sim_\R K_E+B_E$. For the second 
hypothesis, it suffices to show the following equality:
$$
\cO_W=g_*\cO_E(\lceil -(B_E^{<1}) \rceil)
$$
We have a natural inclusion 
$j: \cO_W \to g_*\cO_E(\lceil -(B_E^{<1})\rceil)$ which is
an isomorphism in the generic point of $W$. Since $\cO_W$ 
is reflexive and $g_*\cO_E(\lceil -(B_E^{<1})\rceil)$ is 
torsion free, it is enough to check surjectivity in 
codimenision one points of $W$ (cf. ~\cite[2.iv]{c3flds}). 
For this, we may assume that $W$ is a curve and $X$ is 
a germ at a closed point $P \in \nu(W)$.
If $\lceil -B' \rceil$ is effective, then $\nu(W)$ is 
normal at $P$ and the desired equality holds. If
$\lceil -B' \rceil$ is not effective, then 
$f_*\cO_Y(\lceil -B' \rceil) \subseteq m_{P,X}$.
On the other hand, $R^1 \mu_* \cO_Y(\lceil -B_Y \rceil)$
is torsion free by ~\ref{tf}.(i). Therefore we have
a surjection
$$
\mu_*\cO_Y(\lceil -B' \rceil) \to 
g_*\cO_E(\lceil -B_E \rceil) \to 0.
$$
In particular, $g_*\cO_E(\lceil -B_E \rceil) \subset m_{Q,W}$
for every point $Q \in \nu^{-1}(P)$. This implies that
$\lceil -(B_E^{<1}) \rceil$ contains none of the fibers
$g^{-1}(Q)$ in its support. Consequently, 
$\cO_W=g_*\cO_E(\lceil -(B_E^{<1}) \rceil)$ at $P$.

By construction, $\nu(LCS(W,\omega))$ is contained in the 
union of $(X,B)_{-\infty}$ and all lc centers of $(X,B)$ 
different than $\nu(W)$ (this is the subscheme of $X$ with 
ideal sheaf $\mu_*\cO_Y(\lceil -B' \rceil)$).

\item[(ii)] We may assume $g$ factors as 
$g=\sigma \circ h$, where $\sigma:W' \to W$ is 
a resolution such that $(E,P)\stackrel{h}{\to} (W',Q)\to S$ 
satisfies the assumptions in ~\ref{jda}, 
$B_E$ is supported by $P$, $\Supp(B_E^h)$ has relative 
normal crossings over $W' \setminus Q$ and $h(\Supp(B_E^v)) 
\subseteq Q$.

Define $B_{W'}=\sum b_i Q_i$ by the formulas
$1-b_i=\min_{P_j/Q_i}\frac{1-b_j}{m_{P_j/Q_i}}$, and 
let $M$ be an $\R$-divisor on $W'$ such that
$$
K_E+B_E \sim_\R h^*(K_{W'}+B_{W'}+M)
$$
Since $g_*\cO_E(\lceil -B_E\rceil) \subset \cO_W$, the
negative part of $B_{W'}$ is exceptional over $W$.
Also, $LCS(W',B_{W'}) \subset \sigma^{-1}(LCS(W,\omega))$:
if $b_i \ge 1$, there exists $P_j/Q_i$ such that $b_j \ge 1$,
hence $\sigma(Q_i)=g(P_j) \subset LCS(W,\omega)$.
Note that $B_{div}=\sigma_*B_{W'}$ is the divisorial 
part of adjunction induced by $(X,B)$ on $W$ 
(cf. ~\cite[\S 3]{thesis}).

Since $D=B_E-h^*B_{W'}$ satisfies the hypothesis of 
~\ref{jda}, $M$ is nef/$S$.
In particular, $M+\sigma^*H$ is nef and big/$S$, so there
exists an effective $\R$-divisor $\Delta$ with 
arbitrary small coefficients such that
$M+\sigma^*H \sim_\R \Delta$. We set 
$B_W=\sigma_*(B_{W'}+\Delta) =B_{div}+\sigma_*\Delta$. 
Then $\sigma:(W',B_{W'}+\Delta) \to (W,B_W)$ is 
a crepant birational contraction, hence the claim. 

\end{itemize}
\end{proof}

\begin{thm}\label{jda} ~\cite[Theorem 1]{adjkaw} Let 
$h:(Y,P) \to (X,Q)$ be a projective contraction of 
non-singular varieties endowed with
simple normal crossings boundaries, $Q=\sum Q_i$ and 
$P=\sum P_j$, such that $h^{-1}(Q) \subset P$ 
and $h$ is smooth over $X \setminus Q$. Assume
$X/S$ is a projective morphism, and $D$ is an 
$\R$-divisor on $Y$ with the following properties:
\begin{itemize}
\item[(0)] There exists a non-singular projective
variety $\bar{X}$ endowed with a simple normal 
crossings divisor $\bar{Q}$ such that 
$\bar{X}\setminus X$ is a simple normal crossings
divisor in $\bar{X}$ having simple normal crossings
with $\bar{Q}$, and $Q=\bar{Q}\cap X$.

\item[(1)] $D=\sum d_i Q_i$ is supported by $P$, 
and if we decompose $D$ into horizontal and vertical 
components $D=D^h+D^v$, then $h(\Supp(D^v))\subset Q$
and $\Supp(D^h)$ has relative normal crossings over 
$X\setminus Q$.

\item[(2)] For each $j$, $d_i \le 1-\mult_{P_i}h^*Q_j$ 
if $h(P_i)=Q_j$, and equality holds for some $i$.

\item[(3)] $\lceil -D^h \rceil$ is effective and 
$\cO_{X,\eta_X}\simeq (h_*\cO_X(\lceil -D \rceil))_{\eta_X}$.

\item[(4)] $K_Y+D\sim_\R h^*(K_X+M)$ for some 
$\R$-divisor $M$ on $X$.
\end{itemize}

Then $M$ is nef/$S$.
\end{thm}

\begin{proof} We show that $D$ is a $\Q$-divisor.
Since $K_Y+D \equiv 0$ over the generic point of $X$,
$D^h$ is rational. Over the generic point of each $Q_j$,
$K_Y+D$ is relatively numerically trivial and at least
one component of $D$ has rational coefficients by (2).
The fibers of $h$ are connected hence $D$ is rational 
over the generic points of $Q$.
In particular, $K_Y+D$ is a rational divisor over a big 
open subset of $X$, and is $\R$-linearly equivalent to 
a pull back from $X$. Therefore $D$ is rational 
(cf. ~\cite[3.25]{PLflips}).

We may choose a $\Q$-divisor $M'$ such that $M'\sim_\R M$
and $K_Y+D\sim_\Q h^*(K_X+M)$. Then the statement is just a 
non-compact version of ~\cite[Theorem 1]{adjkaw}. The same 
argument works, since the semi-positivity is a local analytic 
statement (the covering trick holds by our assumption (0)).
Note that ~\cite[Theorem 1]{adjkaw} is stated under the extra
assumption $\lceil -D \rceil\ge 0$, which is however not 
used during the proof (cf. ~\cite[3.5]{thesis}).
\end{proof}


\section{The cone theorem}


We follow the arguments of ~\cite[2-4]{KMM} and
~\cite{KoM}, which we also refer to for references. 

\begin{thm} [Base Point Free Theorem] \label{bpf}
Assume $X/S$ is a projective quasi-log variety. Let 
$L$ be a $\pi$-nef Cartier divisor on $X$ such that:
\begin{itemize}
   \item[(i)] $qL-\omega$ is a $\pi$-ample for 
   some $q\in \R$.
   \item[(ii)] $\cO_{X_{-\infty}}(mL)$ is 
   $\pi|_{X_{-\infty}}$-generated for $m\gg 0$.
\end{itemize}
Then $\cO_X(mL)$ is $\pi$-generated for $m \gg 0$. 
\end{thm}

\begin{proof} We may shrink $S$ to an affine open subset 
without further notice.

1. $\cO_X(mL)$ is $\pi$-generated on $LCS(X)$ for $m\gg 0$.
Set $X'=LCS(X)$. 
The vanishing $R^1\pi_*\cI_{X'}\otimes \cO_X(mL)=0 \ 
(m\ge q)$
implies the surjectivity of the top horizontal map in 
the diagram below:
 \[ \xymatrix{
   \pi^*\pi_*\cO_X(mL)  \ar[r] \ar[d]_{\alpha} & 
   \pi^*\pi_* \cO_{X'}(mL) \ar[d]^{\alpha'}   \\
    \cO_X(mL) \ar[r] &  \cO_{X'}(mL) 
        } \]
If $X'=X_{-\infty}$, $\alpha'$ is surjective for $m\gg 0$ 
by assumption. If $X'\ne X_{-\infty}$, then $X'$ is a 
quasi-log variety, hence $\alpha'$ is surjective 
for $m\gg 0$ by induction. Therefore $\alpha$ is surjective 
on $X'$ for $m\gg 0$.

2. $\cO_X(mL)$ is $\pi$-generated on a non-empty set 
for $m\gg 0$. According to step (1), we may assume 
$LCS(X)=\emptyset$. In particular, $X$ is normal.
\begin{itemize}

\item[(a)] Assume $L$ is $\pi$-numerically trivial. 
Vanishing implies 
that $\pi_*\cO_X(L)$ and $\pi_*\cO_X(-L)$ are non-zero 
\cite{nonvanishing}. Therefore $L$ is trivial, hence 
$\pi$-generated.
	
\item[(b)] Assume $L$ is not $\pi$-numerically trivial. 
Denote $H=qL-\omega$. Using a quasi-log resolution 
of $X$, we can find an $\R$-divisor $D$ on $X$ such that 
$D\sim_\Q c(H+mL), \ 0<c<1$, and $(X,\omega+D)$ has 
qlog canonical singularities, with non-empty LCS 
locus ~\cite{nonvanishing}. Setting $q'=q+cm$, we are 
reduced to Step 1.
\end{itemize}

3. Assume $\cO_X(mL)$ is $\pi$-generated on a non-empty subset 
   containing $LCS(X)$, and denote by $Bsl_\pi|mL|$ the locus 
   $X$ where $\cO_X(mL)$ is not $\pi$-generated. Then 
   $Bsl_\pi|mL|$ 
   is not contained in $Bsl_\pi|m'L|$ for $m'\gg 0$.

Let $f:(Y,B)\to X$ be a quasi-log resolution. For 
$D\in |mL|$ general, we may assume that $f^*D=F+M$ has 
multicrossings support with respect to $(Y,B_Y)$, where $F$ 
is the $\pi$-fixed part and $M$ is reduced. 
Let $c$ be maximal such that $B_Y'=B_Y+cf^*D$ is a 
sub-boundary above $X\setminus X_{-\infty}$. Then 
$f:(Y,B_Y')\to (X,\omega')$ is a quasi-log resolution of a 
quasi-log variety, with $\omega'=\omega+cD$ and 
$X'_{-\infty}=X_{-\infty}$. Moreover, $(X,\omega')$ has a 
qlc center $C$ included in $Bsl_\pi|mL|$. Applying Step 1 with 
$q'=q+cm$, we infer that $\cO_X(m'L)$ is $\pi$-generated on 
$C$ for $m' \gg 0$.

4. The above steps imply that $\cO_X(aL)$ and $\cO_X(bL)$ 
are $\pi$-generated if $a$ and $b$ are very high powers of two 
prime numbers. Since $a$ and $b$ are relatively prime, they 
generate the semigroup $\Z_{\ge N}$ for some $N$. 
Therefore $\cO_X(mL)$ is $\pi$-generated for $m\ge N$.
\end{proof}

\begin{defn} Let $(X/S,\omega)$ be a quasi-log variety, with 
non qlog canonical locus $X_{-\infty}$. Set
$$
\overline{NE}(X/S)_{-\infty}:=
Im(\overline{NE}(X_{-\infty}/S) \to \overline{NE}(X/S)).
$$
For $D\in \Div(X)_\R$, set 
$D_{\ge 0}:=\{z \in N_1(X/S); D\cdot z \ge 0\}$ (similarly
for $>0,\le 0, <0$) and 
$D^\perp:=\{z \in N_1(X/S); D\cdot z = 0\}$. We also use the
notation
$$
\overline{NE}(X/S)_{D\ge 0}:= \overline{NE}(X/S) \cap D_{\ge 0}
$$
and similarly for $>0,\le 0, <0$.
\end{defn}

\begin{defn} An {\em extremal face} of $\overline{NE}(X/S)$ is a 
non-zero subcone $F\subseteq \overline{NE}(X/S)$ such that
$z,z' \in \overline{NE}(X/S), z+z' \in \overline{NE}(X/S)$ imply
that $z,z'\in F$. Equivalently, 
$F=\overline{NE}(X/S) \cap H^\perp$ for some $\pi$-nef
$\R$-divisor $H \in \Div(X)_\R$ 
(called {\em supporting function of} $F$).
An {\em extremal ray} is a $1$-dimensional extremal face.

\begin{itemize}

\item[(i)] An extremal face $F$ is called {\em $\omega$-negative}
if $F\cap \overline{NE}(X/S)_{\omega \ge 0}=\{0\}$.

\item[(ii)] An extremal face $F$ is called 
{\em relatively ample at infinity} if 
$
F\cap \overline{NE}(X/S)_{-\infty}=\{0\}.
$
Equivalently, $H|_{X_{-\infty}}$ is $\pi|_{X_{-\infty}}$-ample
for any supporting function $H\in \Div(X)_\R$ of $F$.

\item[(iii)] An extremal face $F$ is called {\em contractible at 
infinity} if it has a rational supporting function 
$H\in \Div(X)_\Q$ such that 
$H|_{X_{-\infty}}$ is $\pi|_{X_{-\infty}}$-semi-ample.
 
\end{itemize}
\end{defn}

\begin{rem} 
\begin{itemize}
\item[-] Let $F$ be an extremal face which is ample at infinity.
Then $F$ is contractible at infinity if and only if $F$ is rational, 
i.e. it has a supporting function given by a rational divisor.
We will show in the Cone Theorem that if an $\omega$-negative 
extremal face is ample at infinity, then it is contractible at 
infinity.  

\item[-] Any $\omega$-negative extremal face is relatively 
ample at infinity if $\omega$ is relatively nef on 
$X_{-\infty}$ (in particular, if $X_{-\infty}$ is empty).
\end{itemize}
\end{rem}

\begin{defn} Let $F$ be an extremal face of $\overline{NE}(X/S)$.
The {\em contraction of $F$} is a projective morphism
onto a projective variety $Y/S$
\[ \xymatrix{
X \ar[rr]^{\varphi_F} \ar[dr]_{\pi} &  & Y \ar[dl]^{\sigma} \\
         &                    S  &
} \]
satisfying the following properties:
\begin{itemize}
\item[(1)] Let $C$ be an irreducible curve of $X$ such that 
   $\pi(C)$ is a point. Then $\varphi_F(C)$ is a point if
  and only if $[C]\in F$.
\item[(2)] $\cO_Y =({\varphi_F})_* \cO_X$.
\end{itemize}
By Zariski's Main Theorem, such a morphism is unique 
if it exists.
\end{defn}

\begin{thm}[Contraction Theorem] Let $X/S$ be a projective
quasi-log variety. 
Let $F$ be an $\omega$-negative extremal face of 
$\overline{NE}(X/S)$ which is contractible at infinity. 
Then the contraction of the face $F$ exists.
\end{thm}

\begin{proof} Let $H\in \Div(X)$ be a $\pi$-nef divisor
such that $H|_{X_{-\infty}}$ is relatively semi-ample and
$F=\overline{NE}(X/S) \cap H^\perp$. By Kleiman's ampleness
criteria, $aH-\omega$ is $\pi$-ample for some positive 
integer $a$. Scaling $H$, we may assume that its restriction at 
infinity is relatively free. According to the Base Point Free 
Theorem, some multiple of $H$ if relatively free. The Stein 
factorization $\varphi:X/S\to Y/S$ of the associated morphism
satisfies the following properties:
\begin{itemize}
\item[(1)] $H\sim_\Q \varphi^*(A)$ for some relatively ample
 $A\in \Div(Y)_\Q$.
\item[(2)] $\cO_Y=\varphi_*\cO_X$.
\end{itemize}
Since $A$ is relatively ample, it is clear that $\varphi$ 
is the contraction of the face $F$.
\end{proof}

\begin{rem} Let $F$ be an $\omega$-negative extremal face 
which is contractible at infinity. Then $F$ is relatively 
ample at infinity if and only if the associated contraction 
$\varphi_F:X\to Y$ embeds $X_{-\infty}$ into $Y$.
\end{rem}

\begin{lem}\label{polyn} 
Let $P(x,y)$ be a non-trivial polynomial of degree at most $d$, let 
$a$ be a positive integer and let $r$ be either an irrational number, 
or a rational number such that, in reduced form, $ra$ has numerator 
bigger $(d+1)a$.
Then $P(x,y)\ne 0$ for all sufficiently large integral points in the 
strip $\{rax - r< y < ra x\}$.
\end{lem}

\begin{proof} 
If $r$ is not rational, there are integral points of the strip which 
are infinitely close to the line $\{y=rax\}$. If $r$ is rational, let 
$ra=\frac{u}{v}$ be the reduced form decomposition.
The line $\{y=rax-\frac{1}{v}\}$ has infinitely many integral points, 
and it is included in the strip $\{ra x- \frac{r}{d+1}<y<ra x\}$ if 
$u>a(d+1)$.

In both cases, there are infinitely many rays through the origin 
having at least $d+1$ integral points common with the strip 
$\{rax - r< y < ra x\}$. Since $P$ is non-trivial, it cannot vanish 
on more than a finite number of them.
\end{proof}

\begin{thm}[Rationality Theorem] Assume $X/S$ is a 
projective quasi-log variety such that $\omega \in \Div(X)_\Q$.
Let $H$ be a $\pi$-ample Cartier divisor on $X$, and let $r$ 
be a positive number such that
\begin{itemize}
   \item[(i)] $\omega+rH$ is $\pi$-nef, but not $\pi$-ample.
   \item[(ii)] $(\omega+rH) |_{X_{-\infty}}$ is 
        $\pi|_{X_{-\infty}}$-ample.
\end{itemize}
Then $r$ is a rational number, and in reduced form, 
$ra$ has numerator at most $a(\dim X/S +1)$, where 
$a$ is the index of $\omega$.
\end{thm}

\begin{proof} Assume, by contradiction, that $r$ does not
satisfy the required properties. In particular, the strip
$$\cS=\{(x,y)\in \N^2; ra x -r < y< ra x , (x,y) 
\mbox{ large} \}$$ has infinitely many points. 
Set $L(x,y)=xa\omega+yH$. 
The family of Cartier divisors $\{L(x,y)\}_{(x,y)\in \cS}$ 
has the following properties with respect to $(X,\omega)$:

\begin{itemize} 
  \item[(1)] The locus $\Bsl_\pi |L(x,y)|$, where 
$\cO_X(L(x,y))$ is not $\pi$-generated, is independent 
of $(x,y)\in \cS$. We denote this base locus by $\Lambda$.
    \begin{proof} Note first that if $(x,y)$ is a given point 
of $\cS$ and $(kx,ky)$ is a large multiple which does 
not lie in $\cS$, then $L(x',y')-L(kx,ky)$ is $\pi$-ample 
and $\pi$-generated for $(x',y')\in \cS$ large. In particular, 
for  $(x,y)$ given,  $\Bsl_\pi |L(x,y)|$ contains 
$\Bsl_\pi |L(x',y')|$ for $(x',y') \in \cS$ large. 
The claim follows by Noetherian induction.
\end{proof}

\item[(2)] $L(x,y)$ is an adjoint divisor with respect to 
$\omega$ for all $(x,y)$.
\begin{proof}
$L(x,y) -\omega=(xa-1)(\omega+rH)+(y-rax+r)H$ 
is $\pi$-ample for $y>rax-r$. Note that 
$L(x,y)$ is $\pi$-ample for $y>rax$. 
\end{proof}

\item[(3)] $\Lambda \cap (X,\omega)_{-\infty} =\emptyset$ 
and for each qlc center $C$ of $(X,\omega)$, there exists 
$(x,y)$ such that $\cO_C(L(x,y))$ is $\pi|_C$-generated 
on some non-empty subset.
\begin{proof}
Since $L(x,y)$ are adjoint with respect to $\omega$, we can lift
global sections of $\cO_X(L(x,y))$ from $X_{-\infty}$. 
Therefore $\Lambda$ does not intersect the non-qlog canonical 
locus if $\cO_{X_{-\infty}}(L(x,y))$ is relatively generated 
for infinitely many values in $\cS$. 
The line $y=rax$ is relatively ample on $X_{-\infty}$, hence 
Lemma~\ref{polyn} implies the existence of infinitely many 
points $(x,y)$ of $\cS$ for which $L(x,y)|_{X_{-\infty}}$ is 
relatively ample. The same argument as in $(1)$ shows
that $\cO_{X_{-\infty}}(L(x,y))$ are relatively generated 
for large values.
	   
For the latter part, let $C$ be a qlc center of $X$. We may 
assume that $C$ does not intersect $X_{-\infty}$, and $S$ 
is a point. By adjunction, $L(x,y)|_C$ are adjoint, hence 
$$
P(x,y)=\dim H^0(C,\cO_C(L(x,y)))=\chi(C,\cO_C(L(x,y)))
$$ 
is a polynomial of degree at most $\dim C\le \dim X/S$. 
It is non-trivial polynomial, hence $P(x,y)\ne 0$ for 
$(x,y)\in \cS$ by Lemma~\ref{polyn} again.
	         \end{proof}
\end{itemize}
By adjunction, for any family $L(x,y)$ satisfying $(1)-(3)$ 
above, the common base locus $\Lambda$ does not intersect 
$X_{-\infty}$ and does not contain any qlc center of $X$.

If $\Lambda=\emptyset$, then $\cO_X(L(x,y))$ is $\pi$-generated, 
in particular $\pi$-nef. This is a contradiction. Therefore 
$\Lambda$ is non-empty. Let $D$ be a general member of 
$|L(x,y)|$, and choose $0<c \le 1$ maximal such that 
$\omega':=\omega+cD$ has qlog canonical singularities outside 
$X_{-\infty}$. 
Note that $(X,\omega')$ and $(X,\omega)$ have the same 
non-qlog canonical locus, and $(X,\omega')$ has a qlc 
center contained in $\Lambda$.
But $\{L(x,y)\}_{(x,y)\in \cS}$ has the same properties 
$(1)-(3)$ with respect to $(X,\omega')$, hence $\Lambda$ 
cannot contain any qlc center of $(X,\omega')$. 
Contradiction.
\end{proof}

\begin{thm}[Cone Theorem] \label{cone} 
Let $(X/S, \omega)$ be a projective quasi-log variety. 
Let $\{R_j\}$ be the $\omega$-negative extremal rays 
of $\overline{NE}(X/S)$ which are relatively ample at 
infinity. Then
\begin{itemize}
 \item[(i)] 
$
\overline{NE}(X/S)=\overline{NE}(X/S)_{\omega \ge 0}+
\overline{NE}(X/S)_{-\infty}+\sum R_j
$
\item[(ii)] There are only finitely many $R_j$'s included 
in $(\omega+H)_{<0}$, for any relatively 
ample $H \in Div(X)_\R$.
In particular, the $R_j$'s are discrete in the half space 
$\omega_{<0}$.

\item[(iii)] Let $F$ be an $\omega$-negative extremal face 
of $\overline{NE}(X/S)$ which is relatively ample at 
infinity. Then $F$ is a rational face (in particular,
contractible at infinity).
\end{itemize}
\end{thm}

\begin{proof} Assume first that $\omega \in \Div(X)_\Q$.
\begin{itemize}
\item[(1)] If $\dim_\R N_1(X/S)\ge 2$, then 
$$\overline{NE}(X/S)=\overline{NE}(X/S)_{\omega \ge 0}+
\overline{NE}(X/S)_{-\infty} +\overline{\sum_F  F},
$$
where the $F's$ vary among all rational proper $\omega$-negative
extremal faces which are relatively ample at infinity, 
and the overline denotes the closure with respect to 
the real topology.
\begin{proof} Denote by $B$ the right hand side.
If equality does not hold, there exists a separating function
$M\in \Div(X)\setminus \{0\}$, which is not a multiple
of $\omega$ in $N^1(X/S)$, such that $M$
is positive on $B\setminus\{0\}$, but is not relatively nef. 
Since $M$ belongs to the interior of the dual cone of 
$\overline{NE}(X/S)_{\omega \ge 0}$, we can scale it so 
that $M=\omega+H$ for a relatively ample $\Q$-Cartier divisor
$H$.
	 
Let $r > 1$ be the largest real number such that $\omega+rH$ 
is relatively nef, but not ample. In particular, $\omega+rH$ 
is relatively ample on $X_{-\infty}$. By the 
Rationality and Contraction Theorems, $r$ is a rational number 
and the extremal face $F\ne \{0\}$, with supporting function 
$\omega+rH$, can be contracted. If $F$ is proper, it is 
contained in $B$, hence $M$ is relatively ample on $F$. 
This contradicts $r > 1$. Otherwise $\omega+rH$ is trivial
and $M= \frac{r-1}{r} \omega$ in $N^1(X/S)$, which contradicts
the choice of $M$.
\end{proof}

\item[(2)] We may take only proper rays in $(1)$:
\begin{proof} Let $F$ be a rational proper $\omega$-negative
extremal face which is relatively ample at infinity, 
and assume $\dim(F)\ge 2$.
Let $\varphi_F:X\to W$ be the associated contraction, 
so that $-\omega$ is $\varphi_F$-ample. Applying $(1)$ 
to $X/W$ we obtain
$$
F=\overline{NE}(X/W)\setminus \{0\}=
(\overline{NE}(X/W)_{-\infty}+\overline{\sum_G  G}) 
\setminus \{0\},
$$
where the $G$'s are the rational proper $\omega$-negative
extremal faces of $\overline{NE}(X/W)$ which are relatively 
ample at infinity. Since $\varphi_F$ embeds $X_{-\infty}$ into
$W$, $\overline{NE}(X/W)_{-\infty}=0$. The $G$'s are
also $\omega$-negative extremal faces of $\overline{NE}(X/S)$ 
which are contractible at infinity, and $\dim G< \dim F$. 
We obtain by induction
$$
\overline{NE}(X/S)=\overline{NE}(X/S)_{\omega \ge 0}+
\overline{NE}(X/S)_{-\infty}+\overline{\sum  R_j}.
$$
Note that each $R_j$ does not intersect 
$\overline{NE}(X/S)_{-\infty}$.
\end{proof}

\item[(3)] Let $A$ be a relatively ample Cartier divisor 
on $X$. 
Then each $R_j$ is generated by an irreducible reduced 
curve $C_j$, $r_j=\frac{A\cdot C_j}{-\omega \cdot C_j}$ 
is a rational number, and the denominator of $\frac{r_j}{a}$, 
written in reduced form, is at most $a(d+1)$.
Indeed, each $R_j$ is contractible, and the statement follows 
from the Rationality Theorem applied to the contraction 
$\varphi_{R_j}$.

\item[(4)] Let $\{H_i \}_{i=1}^{ \varrho -1}$ be relatively 
ample Cartier divisors on $X$, which together with $\omega$, 
form a basis over $\R$ of $N^1(X/S)$. By $(3)$, 
$R_j \cap \{z; -a\omega \cdot z=1\}$ is included in the 
lattice
$$
\{z; -a\omega \cdot z =1, H_i \cdot z \in (a(d+1)!)^{-1} \Z\}.
$$
Therefore the extremal rays are discrete in the half-space 
$\omega_{<0}$, and the real closure can be omitted. 
We have obtained (i).
	
\item[(5)] We show (ii). Let $H \in \Div(X)_\R$ be relatively
ample. Since
$H-\sum_{i=1}^{ \varrho -1} \epsilon_i H_i$ is ample for 
$0<\epsilon_i\ll 1$, the $R_j$'s included in 
$(\omega +H)_{< 0}$ correspond
to some elements of the above lattice for which
$
\sum_{i=1}^{ \varrho -1} \epsilon_i H_i \cdot z <\frac{1}{a}.
$
They are finite.

\item[(6)] We show (iii). 
The vector space $V=F^\perp \subset N^1(X)$ is defined 
over $\Q$, since $F$ is generated by some of the $R_j$'s. 
There exists a relatively ample divisor $H\in \Div(X)$ 
such that $F \subset (\omega +H)_{<0}$. 
Let $<F>$ be the vector space spanned by $F$, and set 
$$
W_F=\overline{NE}(X/S)_{\omega +H \ge 0} +
\overline{NE}(X/S)_{-\infty}+
\sum_{R_j \not{\subseteq} F} R_j.
$$
Then $W_F$ is a closed cone, $\overline{NE}(X/S)=W_F+F$,
$W_F\cap <F>=\{0\}$, and the supporting functions of $F$
are the elements of $V$ which are positive on 
$W_F\setminus\{0\}$. This is a non-empty open set, thus
contains a rational element which, after scaling, gives a 
relatively nef Cartier divisor $L$ such that 
$F=L^\perp \cap \overline{NE}(X/S)$. Therefore $F$ is 
rational.
\end{itemize}
 
The general case when $\omega \in \Div(X)_\R$ can
be reduced to the rational case via the following trick: 
if $H\in \Div(X)_\R $ is relatively ample and 
$\omega+H \in \Div(X)_\Q$, we can write $H=E+H'$ such 
that $H'\in \Div(X)_\R$ is a relatively ample and 
$(X,\omega':=\omega+E)$ is a quasi-log variety 
with the same qlc centers and non-qlog canonical locus 
as $(X,\omega)$. Therefore $\omega+H=\omega'+H'$, 
$\omega'\in \Div(X)_\Q$ and $(X,\omega)_{-\infty}=
(X,\omega')_{-\infty}$.
In (ii) we may assume that $\omega+H\in \Div(X)_\Q$,
and in (iii) we may replace $\omega$ by
$\omega+H\in \Div(X)_\Q$. As for (i), we have 
 $$
 \overline{NE}(X/S)=\overline{NE}(X/S)_{\omega+H \ge 0}
 +\overline{NE}(X/S)_{-\infty}+\sum_{(\omega+H) \cdot 
 R_j<0} R_j
 $$
 since the same holds for $\omega'+H'=\omega+H$. Letting $H$
 converge to $0$, we obtain (i) using (ii).
\end{proof}

\begin{cor} Let $X/S$ be a projective quasi-log variety 
such that $\omega$ is relatively nef on $X_{-\infty}$. 
If $\omega$ is not relatively nef, there exists an 
$\omega$-negative extremal ray which is relatively ample 
at infinity.
\end{cor}


\section{Quasi-log Fano contractions}


We specialize the results of the previous section to
the equivalent of Fano contractions in our category:

\begin{defn} A {\em quasi-log Fano contraction}
$X/S$ is a relative projective quasi-log variety $X/S$ 
such that $-\omega$ is relatively ample and 
$\cO_S=\pi_*\cO_X$.
\end{defn}

\begin{thm} A projective quasi-log Fano contraction 
$X/S$ has only finitely many $\omega$-negative 
extremal rays $R_j$ which are relatively ample at 
infinity, and  
$
\overline{NE}(X/S)=\overline{NE}(X/S)_{-\infty}+\sum R_j.
$

Furthermore, $NE(X/S)$ is a closed rational polyhedral cone 
spanned by the $R_j$'s, if $X_{-\infty}/S$ has at most
finite fibers.
\end{thm}

\begin{lem}\label{char} 
Assume $X/T\to S/T$ is a diagram of projective 
morphisms such that $X/S$ is a quasi-log Fano contraction.
\begin{itemize}
 \item[(i)] There exists an $\omega$-negative extremal 
face $F$ of $\overline{NE}(X/T)$ which is contractible at 
infinity such that $X/T\to S/T$ is the contraction 
of the face $F$.
 
 \item[(ii)] Let $L\in \Div(X)_K$ such that $L\equiv 0/S$. 
Then there exists $H\in \Div(S)_K$ such that $L\sim_K \pi^*H$, 
if one of the following hold:
\begin{enumerate}
\item[$K=\Z$:] $mL|_{X_{-\infty}}$ is relatively base point 
   free for $m\gg 0$.
\item[$K=\Q$:] $L|_{X_{-\infty}}$ is relatively semi-ample.
\item[$K=\R$:] $X_{-\infty}/S$ has at most finite fibers.
\end{enumerate}

\end{itemize}
\end{lem}

\begin{cor}\label{nsa} 
Let $X/S$ be a quasi-log Fano contraction.
\begin{itemize}
\item[(i)] Assume $L\in \Div(X)_\Q$ is relatively nef, and 
$L|_{X_{-\infty}}$ is relatively semi-ample. Then $L$ is 
relatively semi-ample.
\item[(ii)] Assume $L\in \Div(X)_\R$ is relatively nef, and 
$L|_{X_{-\infty}}$ is relatively ample. Then $L$ is 
relatively semi-ample.
\end{itemize}
\end{cor}

\begin{proof} The first statement follows from the 
Base Point Free Theorem. For (ii), assume $L\in \Div(X)_\R$ 
is relatively nef, and $L|_{X_{-\infty}}$ 
is relatively ample. If $[L]= 0 \in N^1(X/S)$ we just apply 
~\ref{char}.ii.

If $[L]\ne 0 \in N^1(X/S)$, 
$F:=L^\perp\cap \overline{NE}(X/S)$ is a non-trivial face.
By assumption, $F\cap (\overline{NE}(X/S)_{\omega \ge 0}
+ \overline{NE}(X/S)_{-\infty})=\{0\}$. 
Theorem~\ref{cone}.(iii) and the Contraction Theorem imply 
that $F$ is an $\omega$-negative extremal face  
contractible at infinity, and the contraction 
$\varphi_F:X/S\to T/S$ exists. We have $L\equiv 0/T$
and $X_{-\infty}/T$ is an embedding. By ~\ref{char}.ii,
$L\sim_\R \pi^*H$ for some relatively ample 
$H\in \Div(T)_\R$, i.e. $L$ is relatively semi-ample.
\end{proof}

\begin{rem} (cf. Artin's numerical criteria)
Let $\pi:X\to S$ be a projective birational morphism of 
normal varieties, and let $D$ be an effective $\Q$-Cartier
divisor on $X$ such that the following hold:
\begin{itemize}
 \item[-] $(X,B)$ is a log variety.
 \item[-] $-D$ is $\pi$-ample.
 \item[-] For every subscheme $Y\subset X$ supported by
      $\Supp(D)$, any $\pi$-nef Cartier divisor 
      $L\in \Div(Y)$ is $\pi$-semi-ample.
\end{itemize}
Then any $\pi$-nef Cartier divisor $L$ on $X$ is 
$\pi$-semi-ample.
Indeed, $(X/S,B+rD)$ is a quasi-log Fano contraction for 
$r\gg 0$, with non-log canonical locus supported by 
$\Supp(D)$. The claim follows from ~\ref{nsa}(i).
\end{rem}

\begin{thm} \label{u} Let $\pi:X\to S$ be a quasi-log Fano 
contraction, and let $P\in S$ be a closed point. 
\begin{itemize}
  \item[(i)] Assume $X_{-\infty}\cap \pi^{-1}(P)\ne \emptyset$
   and $C$ is a qlc centre such that 
   $C\cap \pi^{-1}(P) \ne \emptyset$ . Then 
	 $C \cap X_{-\infty} \cap \pi^{-1}(P)\ne \emptyset$.
  \item[(ii)] Assume $X$ has qlog canonical singularities. 
     Then the set of all qlc centres intersecting $\pi^{-1}(P)$ 
     has a unique minimal element with respect to inclusion.
\end{itemize}
\end{thm}

\begin{proof} Let $C$ is a qlc center of $X$ such that
$P\in \pi(C)\cap \pi(X_{-\infty})$. By 
Theorem~\ref{adj_van} (with $L=0$),
$X':=C\cup X_{-\infty}$ is a quasi-log variety and
the restriction map $\pi_*\cO_X\to \pi_*\cO_{X'}$ is 
surjective. Since $\cO_S=\pi_*\cO_X$, $X_{-\infty}$ and 
$C$ intersect over a neighborhood of $P$.

Assume now that $X_{-\infty}=\emptyset$, and let $C_1,C_2$ 
be two qlc centers of $X$ such that $P\in \pi(C_1)\cap \pi(C_2)$. 
The union $X'=C_1\cup C_2$ is a quasi-log variety, and 
the same argument implies the surjectivity of the
restriction map $\pi_*\cO_X\to \pi_*\cO_{X'}$. Therefore
$C_1$ and $C_2$ intersect over $P$. 
Furthermore, the intersection $C_1\cap C_2$ is a union of qlc 
centres by Proposition~\ref{conn}. By induction, there exists a 
unique qlc centre $C_P$ over a neighborhood of $P$ such that 
$C_P\subseteq C$ for every qlc centre $C$ with $P\in \pi(C)$.
\end{proof} 


\section{The log big case}


For certain applications, we need to weaken the
projectivity assumption in the Base Point Free Theorem. 

\begin{defn} [M. Reid] Let $X/S$ be a proper quasi-log
 variety. A relatively nef $\R$-Cartier divisor $H$ on $X$ is 
called {\em log big} if $H|_C$ is relatively big for every 
qlc centre $C$ of $X$.
\end{defn}

\begin{thm} [cf. ~\cite{Fk,Fj}] \label{logbig}
Let $X/S$ be a proper quasi-log variety, and 
let $L$ be a relatively nef Cartier divisor
on $X$ with the following properties:
\begin{itemize}
   \item[(i)] $qL-\omega$ is relatively nef and 
   log big for some $q\in \R$.
   \item[(ii)] $\cO_{X_{-\infty}}(mL)$ is 
   relatively generated for $m\gg 0$.
\end{itemize}
Then $\cO_X(mL)$ is $\pi$-generated for $m \gg 0$. 
\end{thm}

The proof is parallel to Theorem~\ref{bpf}. We just need the 
appropriate equivalent of Theorem~\ref{adj_van}:

\begin{thm} Let $X/S$ be a proper quasi-log variety, and 
let $X'$ be the union of $X_{-\infty}$ with a 
union of some qlc centers of $X$. Let $L$ be a Cartier divisor 
on $X$ such that $L-\omega$ is relatively nef and log big. 
Then $\cI_{X'}\otimes \cO_X(L)$ is $\pi_*$-acyclic.    
\end{thm}

This is a formal consequence of the log big extension
of Theorem~\ref{tf}, which we prove below by reduction
to the ample case:

\begin{thm}
Let $f:(Y,B)\to X$ be a proper morphism from an
embedded normal crossings pair, such that $B$ is 
a boundary. Let $L\in Div(Y)$, let 
$\pi:X\to S$ be a proper morphism, and assume 
that $L\sim_\R K+B+ f^*H$ for a nef and log 
big/$S$ $\R$-Cartier divisor $H$ on $X$. Then:
\begin{itemize}
   \item[(i)] 
Every non-zero local section of $R^qf_*\cO_Y(L)$ 
contains in its support the $f$-image of some 
strata of $(Y,B)$.

\item[(ii)] $R^qf_*\cO_Y(L)$ is $\pi_*$-acyclic.
\end{itemize}
\end{thm}

\begin{proof} (1) Assume first that each strata of
$(Y,B)$ dominates some irreducible component of $X$.
Taking the Stein factorization, we may assume that 
$f$ has connected fibers. Assume then that $X$ is 
connected, which implies that $X$ is irreducible and
each strata of $(Y,B)$ dominates $X$.
By Chow's lemma, there exists a proper birational
morphism $\mu:X'/S\to X/S$ such that $X'/S$ is 
projective. Replacing $Y$ by some blow-up, we may
assume that $f$ factors through $\mu$: $f=\mu\circ g$.
Set $\cF=R^qg_*\cO_Y(L)$. Since $\mu^*H$ is nef and big 
over $S$, and $X'/S$ is projective, we may write 
$\mu^*H=E+A$, where $E$ is an effective $\R$-divisor 
such that $B+g^*E$ has multicrossings support and 
$\lfloor B \rfloor=\lfloor B+g^*E \rfloor$, and 
$A\in \Div(X')$ is ample over $S$. From the ample case, 
we infer that $\cF$ is 
$\mu_*$- and $(\pi\circ \mu)_*$-acyclic, and satisfies
(i). Therefore $R^qf_*\cO_Y(L) \simeq \mu_* \cF$ satisfies 
(i) and (ii).
 
(2) We treat the general case by induction on $\dim X$.
We may assume that $Y=Y'\cup Y''$ is a decomposition of 
$Y$ such that $Y'$ is the union of all strata of $(Y,B)$ 
which are not mapped to irreducible components of $X$.
Since $f:(Y'',B'')\to X$ and $L''$ satisfy the assumption 
in (1), the long exact sequence of 
$0\to j_*\cO_{Y''}(L'') \to \cO_Y(L)
\to \cO_{Y'}(L) \to 0$ with respect to $f_*$, 
breaks up into short exact sequences
$$
0\to R^q f_*\cO_{Y''}(L'') \to R^q f_*\cO_Y(L) \to 
R^qf_*\cO_{Y'}(L)\to 0.
$$
Since (i) and (ii) hold for the first and third member
by case (1) and by induction on dimension respectively, 
they hold for $R^q f_*\cO_Y(L)$ also.
\end{proof}


\end{document}